  \documentclass[11pt]{article}
  \usepackage[top=1in,bottom=1in,left=1.5in,right=1.5in]{geometry}
   \usepackage{amsmath,amssymb}
   \usepackage{latexsym}
   \usepackage[dvips]{graphicx}
   
 \newcommand{\submodularf}{w}
   
 \newcommand{\interval}{\mathcal{E}}
   \newcommand{\C}{\mathbb{C}}
   
   \newcommand{\N}{\mathbb{N}}
   
   \newcommand{\R}{\mathbb{R}}

   \newcommand{\U}{\mathbf{U}}
   \newcommand{\uu}{\mathbf{u}}
   
   \newcommand{\V}{\mathbf{V}}
   
   \newcommand{\W}{\mathbf{W}}
   \newcommand{\x}{\mathbf{x}}
   
   \newcommand{\y}{\mathbf{y}}

   \newcommand{\0}{\mathbf{0}}
   \newcommand{\1}{\mathbf{1}}

   \newcommand{\cE}{\mathcal{E}}

   \newcommand{\cH}{\mathcal{H}}

   \newcommand{\cU}{\mathcal{U}}

   \newcommand{\cX}{\mathcal{X}}

   \newcommand{\rH}{\mathrm{H}}

   \newcommand{\lan}{\langle}
   \newcommand{\ran}{\rangle}
   \newcommand{\an}[1]{\lan#1\ran}
   \def\diag{\mathop{{\rm diag}}\nolimits}
   \newcommand{\hs}{\hspace*{\parindent}}
   \newcommand{\proof}{\hs \textbf{Proof.\ }}
   
   \newcommand{\tr}{\mathop{\mathrm{tr}}\nolimits}

   \newcommand{\qed}{\hspace*{\fill} $\Box$\\}

   \newcommand{\Id}{\mathrm{Id}}

   \newcommand{\rC}{\mathrm{C}}

   \newcommand{\rS}{\mathrm{S}}

   \newcommand{\rank}{\mathrm{rank\;}}

   \newtheorem{theo}{\bfseries \hs Theorem}
   
   \newtheorem{prop}[theo]{\bfseries \hs Proposition}
   
   \newtheorem{lemma}[theo]{\bfseries \hs Lemma}
   \newtheorem{corol}[theo]{\bfseries \hs Corollary}
   
   \newtheorem{algo}[theo]{\bfseries \hs Algorithm}
   \newtheorem{example}[theo]{\bfseries \hs Example}
   \newtheorem{rem}[theo]{\bfseries \hs Remark}
   
   \numberwithin{equation}{section} 

  \renewcommand{\span}{\mathrm{span}}

\renewcommand{\geq}{\geqslant}
\renewcommand{\leq}{\leqslant}
\renewcommand{\ge}{\geqslant}
\renewcommand{\le}{\leqslant}

\begin{document}

  \title{Submodular spectral functions of principal submatrices\\of a hermitian matrix, extensions
  and applications}
  \author{
  S. Friedland\footnotemark[1], and\;
  S. Gaubert\footnotemark[2]
  }
  \renewcommand{\thefootnote}{\fnsymbol{footnote}}
  \footnotetext[1]{
  Department of Mathematics, Statistics and Computer Science,
  University of Illinois at Chicago, Chicago, Illinois 60607-7045,
  USA, \texttt{friedlan@uic.edu},
  }
  \footnotetext[2]{INRIA and CMAP, \'Ecole Polytechnique, 91128 Palaiseau C\'edex, \texttt{stephane.gaubert@inria.fr} }
  \renewcommand{\thefootnote}{\arabic{footnote}}
  \date{June 17, 2012}
  \maketitle
  \begin{abstract}
  We extend the multiplicative submodularity of the principal determinants
  of a nonnegative definite hermitian matrix to other spectral functions.
  We show that if $f$ is the primitive of a function that is operator
  monotone on
  an interval containing the spectrum of a hermitian matrix $A$, then
  the function $I\mapsto {\rm tr} f(A[I])$ is supermodular, meaning
  that ${\rm tr} f(A[I])+{\rm tr} f(A[J])\leq {\rm tr} f(A[I\cup J])+{\rm
  tr} f(A[
  I\cap J])$, where $A[I]$ denotes the $I\times I$ principal submatrix of $A$.
  We discuss extensions to self-adjoint operators on infinite dimensional Hilbert space and
  to $M$-matrices.  We discuss an application to CUR approximation of nonnegative hermitian matrices.
  \end{abstract}



  \section{Introduction and statement of the main result}\label{intro}

  Let $m$ be a positive integer and denote $[m]:=\{1,\ldots,m\}$.
 A real valued function $\submodularf:2^{[m]}\to \R$ defined on all subsets of $[m]$ is called \emph{nondecreasing} if $\submodularf(I)\le \submodularf(J)$ when $I\subset J\subset [m]$.
 It is \emph{submodular} if
 \[ \submodularf(I)+\submodularf(J)\ge \submodularf(I\cup J)+\submodularf(I\cap J)
 \]
 for any two subsets $I,J$ of $[m]$.
 The function $\submodularf$ is called \emph{nonincreasing} or~\emph{supermodular}
 whenever $-\submodularf$ is nondecreasing or submodular, respectively. A function that is both submodular and supermodular is called~\emph{modular}.

 The importance of submodular functions in combinatorial optimization is well known. Several polynomial time algorithms to minimize a submodular function
 under a matroid constraint are known, we refer the reader to the
 survey~\cite{iwata} for more information. The maximization
 of a submodular function under a matroid constraint, and specially,
 under a cardinality constraint,
 $\nu_k(\submodularf):=\max_{I\subset [m], |I|\le k} \submodularf(I)$,
 is also of great interest. 
  For some submodular functions $\submodularf$ the latter problem is NP-hard.
 However, a classical result~\cite{NWF78} shows
 that when $\submodularf$ is nondecreasing and submodular, the greedy algorithm
 allows one to compute an approximation $\nu_k^G(\submodularf)$ of $\nu_k(\submodularf)$ which is
 such that $\nu_k^G(\submodularf)\ge (1-e^{-1})\nu_k(\submodularf)$.
 See~\cite{vondrak} for recent developments regarding submodular maximization.

  Denote respectively by $\C^{m\times m}\supset\rH_m\supset \rH_m(\interval)$ the space of $m\times m$ complex valued matrices, the space of $m\times m$ hermitian matrices, and the subset of $A\in \rH_m$ whose eigenvalues lie in the interval $\interval\subset \R$.
  For $I\subset [m]$ denote by $A[I]$ the principal submatrix of $A$,
  obtained from $A$ by deleting the rows and columns in the set $[m]\setminus I$.
  For $A,B\in \rH_{m}$ denote by $\succeq$ the Loewner ordering, so that $A\succeq B$ if  $A-B\in \rH_{m}([0,\infty))$. We also write
  $A\succneqq B$ or $A\succ B$ if $A-B
  \in\rH_{m}([0,\infty))\setminus \{0\}$ or $A-B\in \rH_{m}((0,\infty))$ respectively.

  Recall that the principal minors of a nonnegative definite matrix satisfy the multiplicative submodularity property:
  \begin{equation}\label{multsubprop}
  \det A[I\cup J]\det A[I\cap J]\le \det A[I]\det A[J], \textrm{ where } I,J\subset [m],\;
  A\in\rH_m([0,\infty)).
  \end{equation}
  In other words, the function $\log(\cdot,A): 2^{[m]}\to \R$ given by
  \begin{equation}\label{deffunc}
  \log (I,A):=\log\det A[I], \quad I\subset [m],\; A\in\rH_{m}((0,\infty))
  \end{equation}
 is submodular.
 This inequality has arisen in the work of several authors.
 It goes back to Gantmacher and Kre{\u\i}n~\cite{gantmacherkrein}
 and Kotelyanski{\u\i}~\cite{kotelyanskii}, see the discussion by Ky Fan~\cite{KyFan67,KyFan68}. It
 can also be found in~\cite{KK83,johnsonbarret85}.
 The classical Hadamard-Fischer inequality for the principal minors of nonnegative definite matrices is obtained when $I\cap J=\emptyset$, understanding
 that $\det A[\emptyset]=1$. We refer the reader to~\cite{fallatjohnson99}
 for a survey of determinantal inequalities, and to~\cite{coverthomas} for their
 relation with information theory.
 It is well known that the inequality \eqref{multsubprop} hold also for $M$-matrices, e.g. \cite{Car67}.

 In this paper, we derive a general submodularity result for spectral functions of hermitian matrices, for some $p$-trace of $M$-matrices and some extensions to self adjoint positive operators on infinite dimensional separable Hilbert spaces.

 We first summarize 
our results on hermitian matrices.
 Recall that if $A=UDU^*$, where $U$ is unitary and $D=\diag(\lambda_1,\dots,\lambda_m)$ is the diagonal matrix with the eigenvalues of $A$ on the diagonal,
 the matrix  $f(A)$ is defined to be $U\diag(f(\lambda_1),\ldots,f(\lambda_m))U^*$. Note in particular  that $\tr f(A)=f(\lambda_1)+\dots+f(\lambda_m)$.

 Recall that a real function $f$
 is \emph{operator monotone}
 on the interval $\interval\subset \R$
 if for all $m\geq 1$ and for all $A,B\in\rH_m(\interval)$,
 \[
 A\preceq B\implies f(A) \preceq f(B) \enspace  .
 \]
 An \emph{operator convex} function on the interval $\interval$ is defined
 by requiring that
 \[ tf(A)+(1-t)f(B)\succeq f(tA+(1-t)B)
 \]
 for each $t\in [0,1]$.
 A function $f$
 is \emph{operator antitone} (resp.\ \emph{operator concave})
 on an interval
 if $-f$ is {operator monotone} (resp.\ {operator convex})
 on the same interval. Operator monotone and operator convex
 functions are characterized by Loewner theory~\cite{loewner}.
In particular,
 integral representations like the one of~\eqref{e-integral} below are known.
 We refer the reader to~\cite{Bha97} for more background.
 Matrix inequalities involving operator
 monotone functions can be found in~\cite{andozhan,zhan}.

 The main result of this paper is the following.
 \begin{theo}\label{theo-newmain}
 Let $f$ be a real continuous function defined on an interval $\interval$ of $\R$,
 and assume that $f$ is the primitive of a function that is operator
 monotone on the interior of $\interval$. Then, for every
 $m\times m$ hermitian matrix $A$ with spectrum in $\interval$,
 the function $2^{[m]}\to \R$
 \[
 I\mapsto \tr f(A[I])
 \]
 is supermodular.
 \end{theo}
 It is understood that $\tr f(A[\emptyset]):=0$ for every function $f$.

 Recall that the primitive of an operator monotone function is
 operator convex on the same interval, but that not all operator convex
 functions are obtained in this way. We shall
 see that the conclusion of this theorem no longer holds if $f$ is only assumed
 to be operator convex.

 Let us mention some immediate applications.
 The derivative of the map $t^p$, namely, $pt^{p-1}$,
 is known to be operator antitone on $(0,\infty)$ for $0<p\leq 1$, and
 operator monotone on $[0,\infty)$ for $1\leq p\leq 2$.
 The derivative of the map $t\log t$, namely, $1+\log t$,
 is also known to be operator monotone on $(0,\infty)$.
 (See~\cite{Bha97}.)
 Hence, the next result readily follows from the main theorem.
 \begin{corol}\label{coro-powerentrop}
 Let $A$ be a $m\times m$ nonnegative definite hermitian matrix.
 Then, for all $I,J\subset [m]$,
 \begin{eqnarray}\label{powerentro01}
 \tr A[I]^p+\tr A[J]^p &\geq& \tr A[I\cup J]^p+\tr A[I\cap J]^p,
 \quad \text{for}\quad 0\leq p \leq 1 \enspace ,\\
 \label{powerentro12}
 \tr A[I]^p+\tr A[J]^p &\leq& \tr A[I\cup J]^p+\tr A[I\cap J]^p,
 \quad \text{for}\quad 1\leq p \leq 2 \enspace ,
 \end{eqnarray}
 and
 \begin{eqnarray}\label{xlogxsub}
 \lefteqn{\tr \big(A[I]\log A[I]\big)+\tr \big(A[J]\log A[J]\big) \leq}\\
 \notag
 &&\qquad\qquad\qquad\tr \big(A[I\cup J]\log A[I\cup J]\big)+\tr \big(A[I\cap J]\log A[I\cap J]\big)\enspace .
 \end{eqnarray}\qed
 \end{corol}
 We write $A[I]^p$ for $(A[I])^p$.
 It is understood that $A^0:=\lim_{p\searrow 0}A^p$, for every nonnegative definite matrix $A$. In particular, $A^0$ is the identity matrix $\Id$
 if $A$ is positive definite, and, if $A$ is nonnegative definite, $\tr A^0=\operatorname{rank} A$.

 We now survey briefly the contents of our paper.
 In \S2 we prove Theorem \ref{theo-newmain}.  In \S3 we discuss $M$-matrices.  We show the inequalities \eqref{multsubprop}, \eqref{powerentro01} and \eqref{xlogxsub}.  Moreover $\tr A[I]^p$ is a supermodular function for each $p<0$.
 In \S4 we discuss the extensions of Theorem~\ref{theo-newmain} to
the space $\rS(\cH)$ of self adjoint operators on a separable Hilbert space $\cH$.
We consider the lattice of closed subspaces $\cU$ in $\cH$, with
joint and meet operations $\operatorname{clo}(\U+\V)$ and $\U\cap\V$,
where $\operatorname{clo}$ denotes the closure of a set,
together with its sublattice of finite dimensional subspaces $\cU_f$.
 Let $P(\U)\in \rS(\cH)$ be the orthogonal projection on $\U\in\cU$.
 Associate with each $\U\in\cU$ the operator $A(\U)$ which is the restriction of $P(\U)AP(\U)$ to $\U$.
Thus, $A(\U)$ is an analog of $B[I]$ for a hermitian matrix $B$.  For $\U\in \cU_f$ it is straightforward to show that $w(\U):=\tr f(A(\U))$
 is submodular under the conditions given in \S2.  For an infinite dimensional subspace $\U$, we restrict out discussion to $w_p(\U):=\tr(A(\U))^p$
 and the von Neumann entropy $\tilde w(\U):=-\tr(A(\U)\log A(\U)$, under the assumption that $A$ is positive and compact.  Assuming that $\tr A^p<\infty$  we show that $w_p(\U)$ is submodular for $p\in (0,1)$ and supermodular for $p\in (1,2)$.  We show that $\tilde w$ is submodular
 it $-\tr (A\log A)<\infty$.  In \S5 we discuss the CUR approximation \cite{GTZ97} of nonnegative definite hermitian matrix.
 The main problem here is to find a good approximation to the maximum of $\tr(\log A[I])$ on all subsets $I$ of $[m]$ of cardinality  $k$.
 We discuss briefly the obvious greedy algorithm for this problem, give a simple condition where $\tr (\log A[I])$ is nondecreasing,
 and give an estimate for the $CUR$ approximation obtained by the greedy algorithm in the general case.
 In \S6 we give examples to show that in general the results of \S2 are best possible.


{\em Note added to the arXiv postprint}.\/
The authors thank D. Petz for having brought to their
attention his work with K. Audenaert and F. Hiai (Strongly subaddtive functions,
 \emph{Acta Math. Hungar.} 128(4):386--394, 2010), after this paper was published on line in \emph{Linear Algebra Appl.} (December 2011, doi:10.1016/j.laa.2011.11.021).
Theorem~1 of this paper is a slightly more general version of
their Theorem 4.1 (the latter corresponds to the case in which the interval
is $(0,\infty)$). Some results of the present Section~4
concerning the finite dimensional case (Theorem~5)
can also be thought of as an extension of theirs. Their work
is motivated by quantum information theory.

 \section{Proof of Theorem~\ref{theo-newmain}}
 Observe first that if all the eigenvalues of $A$ belong to the interval
 $\interval$, so do the eigenvalues of the principal submatrix $A[I]$. Hence,
 the function $I\mapsto \tr f(A[I])$ is well defined for all $I\subset [m]$.

 It suffices to consider the case in which the interval $\interval$ is bounded.
 Moreover, if the result is established for every matrix $A$ the spectrum
 of which is included in the interior of $\interval$, arguing by density, the
 result must also hold whenever the spectrum of $A$ is included in $\interval$. Hence,
 we may assume that $\interval$ is open. Since
 the property to be established is invariant by a translation and a scaling of
 the interval $\interval$, we finally assume that $\interval=(-1,1)$. Then, a theorem
 of Loewner (Corollary~V.4.5 in~\cite{Bha97}) shows that every operator
 monotone function $g$ on $(-1,1)$ can be written as
 \begin{equation}\label{e-integral}
 g(t) = g(0)+g'(0)\int_{-1}^1 \frac{t}{1-\lambda t} d\mu(\lambda)
 \end{equation}
 where $\mu$ is a probability measure on $[-1,1]$. Moreover, $g'(0)\geq 0$.

 Assume now that $f$ is a primitive of $g$, so that
 \[
 f(t)=a+tg(0)+g'(0)\int_{0}^t ds \int_{-1}^1 \frac{s}{1-\lambda s} d\mu(\lambda)
 \enspace ,
 \]
 for some constant $a$.
 Observe that for $|t|<1$, the denominator $1-\lambda s$ is bounded below by
 $1-|t|$, and so, the above double integral is absolutely convergent.
 Applying Fubini's theorem, we get
 \begin{eqnarray*}
 f(t)&=&
 a+tg(0)+g'(0)\int_{-1}^1 d\mu(\lambda)\int_{0}^t \frac{s}{1-\lambda s} ds\\
 &=&a+tg(0)+g'(0)\int_{-1}^1 \varphi(\lambda,t) d\mu(\lambda)
 \end{eqnarray*}
 where
 \[
 \varphi(\lambda,t):= -\frac{t}{\lambda}-\frac{1}{\lambda^2}\log(1-\lambda t)
 \enspace .
 \]
 Although $\lambda$ appears at the denominator, this expression defines
 a function of $(\lambda,t)$ that extends continuously to $[-1,1]\times (-1,1)$.
 (In particular, $\varphi(0,t)=\frac{1}{2}t^2$.)

 For all $\lambda\in [-1,1]$, consider now the functions  $\bar{w}$ and $w_\lambda$ from $2^{[m]}\to\R$,
 \[
 w_\lambda(I):=\tr \varphi(\lambda,A[I]) \enspace ,
 \qquad \bar{w}(I):=\tr (a\Id+g(0)A)[I] \enspace,
 \]
 so that
 \begin{equation}
 \tr f(A[I]) = \bar{w}(I)+ g'(0)\int_{-1}^1w_\lambda(I)d\mu(\lambda) \enspace .
 \label{e-rep}
 \end{equation}
 Observe first
 that the function $I\mapsto -\frac{1}{\lambda}\tr A[I]$
 is modular, for all $\lambda\neq 0$.
 The eigenvalues of the matrix $\Id-\lambda A$
 belong to the interval $[1-\lambda,1+\lambda]$, which implies
 that this matrix is positive definite.
 Since the sum of supermodular functions is supermodular,
 it follows from the multiplicative submodularity property
 of the determinant, Eqn~\eqref{multsubprop}, that
 \[ w_\lambda(I)=-\frac{1}{\lambda}\tr A[I]-\frac{1}{\lambda^2}\log(I,\Id-\lambda A)
 \]
 is supermodular, as soon as $\lambda\neq 0$. Since $w_\lambda(I)$ depends continuously of $\lambda$, the same is true when $\lambda=0$. Note finally that the map $\bar{w}$ is modular. Since the
 supermodularity property is preserved by taking positive linear combinations
 and integrals with respect to positive measures, the result follows
 from the representation~\eqref{e-rep}.\qed

 \section{Submodularity and super-modularity inequalities for M-matrices}
 Denote by $\R^{m\times m}_+$ the set of nonnegative matrices.  For $A,B\in \R^{m\times m}$ we denote $A\ge B$ if $A-B\in\R_+^{m\times m}$.
 We also write $A\gneqq B$ and $A>B$ if $A-B\in\R^{n\times n}_+\setminus\{0\}$ and $A-B$ has positive entries respectively.
 For $B\in \R^{m\times m}$ denote by $\rho(B)$ the spectral radius of $B$.  The Perron-Frobenius theorem yields that
 $\rho(B)$ is an eigenvalue of a matrix $B\in\R^{m\times m}_+$.
 Recall that $A\in\R^{m\times m}$ is called an {\em $M$-matrix} if all off-diagonal entries of $A$ are nonpositive and all the principal minors
 of $A$ are nonnegative.  Equivalently, $A$ is an $M$-matrix if $A=s \Id - B$ for some $B\in\R_+^{m\times m}$ and $s\ge \rho(B)$.
 In this section we assume that $A$ is an $M$-matrix and $B\in \R^{n\times n}_+$ unless stated otherwise.
 Note that $A=s \Id-B$ is invertible if and only if $s>\rho(B)$.  Assume that $A$ is an invertible $M$-matrix.
 Then for any $p\in\R$ we define
 \[A^p=s^p\sum_{i=0}^{\infty}  {p \choose i}(-s)^{-i}B^i, \qquad \log A=(\log s)\Id -\sum_{i=1}^{\infty} \frac{1}{is^i}B^i.\]
Since $s>\rho(B)$, a standard argument shows that the above series absolutely converges.
When $A$ is singular, we define, for every $p\geq 0$,
 $\tr A^p:=\lim_{t\to s^+}\tr (t\Id -B)^p$. Indeed,
one readily checks that the spectrum of the matrix
$(t\Id -B)^p$ as a pointwise limit as $t\to s^+$,
if $p\geq 0$, so that limit of $\tr (t\Id -B)^p$
as $t\to s^+$ does exist.
(However, we warn the reader of the abusive character of the notation
$\tr A^p$: the matrix $A^p:=\lim_{t\to s^+}(t\Id-B)^p$
may not exist if $0$ is not a semi-simple eigenvalue of $A$.)
 Since $\rho(B[I])\le \rho(B)$ for any $I\subset [n]$, it follows that
 $A[I]$ is also an $M$-matrix.
 Furthermore, if $A$ is invertible then $A[I]$ is invertible.
 We finally agree that if $A$ is singular then $\tr \log A=-\infty$ and $\tr A^p=\infty$ if $p<0$.

 The main result of this section is the following.
 \begin{theo}\label{mmatriceineq}  Let $A$ be an $M$-matrix.
 Then, for all $I,J\subset [m]$,
 \begin{eqnarray*}
 \tr A[I]^p+\tr A[J]^p &\geq& \tr A[I\cup J]^p+\tr A[I\cap J]^p,
 \quad \text{for}\quad 0\leq p \leq 1 \enspace ,\\
 \tr A[I]^p+\tr A[J]^p &\leq& \tr A[I\cup J]^p+\tr A[I\cap J]^p,
 \quad \text{for}\quad p<0 \quad\text{or}\quad 1\leq p\leq 2\enspace ,
 \end{eqnarray*}
 and
 \begin{eqnarray*}
 \lefteqn{\tr \big(A[I]\log A[I]\big)+\tr \big(A[J]\log A[J]\big) \leq}\\
 &&\qquad\qquad\qquad\tr \big(A[I\cup J]\log A[I\cup J]\big)+\tr \big(A[I\cap J]\log A[I\cap J]\big)\enspace .
 \end{eqnarray*}
 \end{theo}
 The proof of this theorem relies on the following lemma.
 %
 \begin{lemma}\label{supmodngtmat}  For any $B\in \R^{m\times m}_+$
 \[\tr B[I]^n +\tr B[J]^n\le \tr B[I\cup J]^n +\tr B[I\cap J]^n\]
 for every integer $n\ge 1$, and the equality holds for $n=1$.
 \end{lemma}
\proof
Consider the complete directed graph with $m$ nodes. Define
a {\em closed walk} of {\em length} $n$
to be a sequence
$\alpha=(i_1,\ldots,i_{n+1})$ of elements of $[m]$ such that $i_{n+1}=i_1$.
We say that $\alpha$ is {\em included} in $I$, and we write $\alpha\subset I$,
if $i_1,\ldots,i_{n+1}\in I$.
The {\em weight} of this walk is $|\alpha|:=B_{i_1i_2}\dots B_{i_{n}i_{n+1}}$.
By a classical result~\cite[Th.~4.7.1]{Sta97}, $\tr B[I]^n$ is the sum of the weights of all closed walks of length $n$ included in $I$.
By a disjunction of cases, we deduce that
\begin{eqnarray}
\tr B[I\cup J]^n& = &
\sum_{\begin{subarray}{l}\alpha \subset I\\ \alpha\not\subset J\end{subarray}} |\alpha|
+
\sum_{\begin{subarray}{l}\alpha \subset  I\\ \alpha\subset J\end{subarray}} |\alpha|
+
\sum_{\begin{subarray}{l}\alpha \subset  J\\ \alpha\not\subset I\end{subarray}} |\alpha|
+
\sum_{\begin{subarray}{l}\alpha \subset  I\cup J\\ \alpha\not\subset I\\ \alpha\not\subset J \end{subarray}} |\alpha| \enspace,\label{e-disjunction}
\end{eqnarray}
with the convention that all sums are restricted to the walks $\alpha$ of length $n$.
Adding the two first sums in~\eqref{e-disjunction} yields $\tr B[I]^n$.
Moreover, adding $\tr B[I\cap J]^n$ to the third sum in~\eqref{e-disjunction} yields $\tr B[J]^n$.
We conclude that
\[
\tr B[I\cup J]^n+\tr B[I\cap J]^n
-\tr B[I]^n-\tr B[J]^n = \sum_{\begin{subarray}{l}\alpha \subset I\cup J\\ \alpha\not\subset I\\ \alpha\not\subset J \end{subarray}} |\alpha|  \geq 0 \enspace,
\]
since the entries of the matrix $B$ are nonnegative. Moreover, when $n=1$,
the latter sum trivially vanishes.
\qed

 \textbf{Proof of Theorem \ref{supmodngtmat}}.
 Suppose first that $s>\rho(B)$.
 We first consider the function $\tr A[I]^p$.
We have
\[
\tr A[I]^p = s^p\tr \Id[I]
- s^{p+1} p \tr B[I]
+ w'_p(I)
\]
where
\[
w'_p(I):=\sum_{i\geq 2} (-1)^is^{p+i} \frac{p(p-1)\cdots(p-i+1)}{i!} \tr B[I]^i
\enspace .
\]
Observe that for $p\in (-\infty,2]$, all the coefficients
of the latter sum have the same sign (or vanish), and this sign
is positive if $p\in (-\infty,0)\cup (1,2)$ and negative
if $p\in (0,1)$. Hence, it follows from
Lemma~\ref{supmodngtmat} that the function $w'_p$ is
supermodular in the former case and submodular in the latter
case. Since $I\mapsto \tr A[I]^p$ is the sum
of the modular function
$I\mapsto s^p\tr \Id[I]
- s^{p+1} p \tr B[I]$ and of $w'_p$, the announced submodularity and supermodularities of the map $I\mapsto \tr A[I]^p$ hold for
this matrix $A$.


Similarly, since $-(1-x)\log(1-x)-x$ has negative Taylor coefficients, it follows as above that $-\tr (A[I]\log A[I])$ is submodular.

 Assume now that $s=\rho(B)$, i.e., that $A$ is a singular $M$-matrix.  Let $A_t=t\Id-B$ for $t>\rho(B)$.
 Letting $t\searrow \rho(B)$ we deduce all the above inequalities for $A$ except for $p=0$.
 Letting $p\searrow 0$ we deduce that $\tr A[I]^0$ is a submodular function.
 (Note that we cannot use $p<0$ since $\tr A^p$ in this case in $\infty$.)
 \qed

 Note that since $\log(1-x)$ has negative Taylor coefficients, it follows as above that the function $I\mapsto \tr (\log A[I])$ is submodular, i.e.,
 we have obtained yet another proof of the inequality~\eqref{multsubprop} for $M$-matrices \cite{Car67}.

 \section{Submodular functions on the lattice of closed subspaces of a Hilbert space}
Let $\cH$ be a separable Hilbert space over $\C$
with the inner product $\an{\x,\y}$ and the norm $\|\x\|=\sqrt{\an{\x,\x}}$.
 Denote by $\rS(\cH)\supset\rS_+(\cH)$ the space of bounded self adjoint linear operators $A:\cH\to\cH$ and the cone of bounded positive self adjoint operators. So $\an{A\x,\y}=\an{\x,A\y}$ for $A\in\rS(\cH)$ and $\an{A\x,\x}\ge 0$ for $A\in\rS_+(\cH)$ for all $\x,\y\in\cH$.
 For $A\in \rS(\cH)$ the operator norm $\|A\|$ is given by $\sup_{\|\x\|\le 1} |\an{A\x,\x}|$.
 More precisely, the spectrum of $A\in \rS(\cH)$ lies in the interval $[a(A),b(A)]$ where
 \begin{equation}\label{specint}
 a(A)=\inf_{\|\x\|\le 1} \an{A\x,\x},\; b(A)=\sup_{\|\x\|\le 1} \an{A\x,\x},\; \|A\|=\max(|a(A)|,|b(A)|).
 \end{equation}

 Assume that $A\in\rS(\cH)$.
 We first discuss the analog of $A[I]$ when $I$ is a finite subset of $\N$ of cardinality $l$.
 Let $\U$ be an $l$-dimensional subspace of $\cH$.  Denote by $P(\U)\in \rS_+(\cH)$ the orthogonal
 projection on $\U$.  So $\cH=\U\oplus \U^{\perp}=\U\oplus(\Id-P(\U))\cH$.  Consider the operator $P(\U)AP(\U)$.
 Clearly $P(\U)AP(\U)\U^{\perp}=\{\0\}, P(\U)AP(\U)\U\subseteq \U$.  Denote by $A(\U)$ the restriction of $P(\U)AP(\U)$
 to $\U$.  As $\|P\x\|\le \|\x\|$ for $\x\in\cH$ the characterization \eqref{specint} yields
 \begin{equation}\label{specintP}
 a(A)\le a(A(\U))\le b(A(\U))\le b(A).
 \end{equation}
In particular, if the spectrum of $A$ lies in a given interval $\cE$,
so does the spectrum of $A(\U)$.

 Let $\V$ be another finite dimensional vector space of $\cH$.  Then $A(\U),A(\V),A(\U+\V),A(\U\cap\V) $ are the analogs
 of $A[I],A[J],A[I\cup J],A[I\cap J]$ respectively.

 Indeed, let $\W=\U+\V$ and choose an orthonormal basis $W:=\{\x_1,\ldots,\x_m\}$ in $\W$, such that $X:=\{\x_1,\ldots,\x_q\}$ is a basis in $\U$,
 $Z:=\{\x_p,\x_{p+1},\ldots,\x_q\}$ is a basis in $\U\cap\V$  and $Y:=\{\x_p,\x_{p+1},\ldots,\x_m\}$ is a basis in $\V$.
 It is easy to see that the operators $A(\W), A(\U),A(\V), A(\U\cap \V)$ are represented in the above bases by the following matrices
 \begin{eqnarray*}
 &&A(W)=[\an{A\x_i,\x_j}]_{i=j=1}^m, \qquad A(X)=[\an{A\x_i,\x_j}]_{i=j=1}^q, \\
 &&A(Y)=[\an{A\x_i,\x_j}]_{i=j=p}^m,  \qquad A(Z)=[\an{A\x_i,\x_j}]_{i=j=p}^q.
 \end{eqnarray*}

 Denote by $\cU\supset \cU_f$ the set of all closed subspaces of $\cH$ and the subset of all
 finite dimensional subspaces, respectively.  A function $w: \cU_f\to \R$
is said to be {\em submodular} if
 $w(\U)+w(\V)\ge w(\U+\V)+w(\U\cap \V)$ for all $\U,\V\in\cU_f$.  It is
{\em supermodular} if $-w$ is submodular, and {\em modular} if it is both submodular and supermodular.
The following result is an immediate consequence of Theorem~\ref{theo-newmain}
and of the previous observations.
 \begin{theo}\label{theo-newmainid}
 Let $f$ be a real continuous function defined on an interval $\interval$ of $\R$,
 and assume that $f$ is the primitive of a function that is operator
 monotone on the interior of $\interval$.
 Let $\cH$ be an infinite dimensional separable Hilbert space space.
 Then, for every $A\in\rS(\cH)$ with spectrum in $\interval$,
 the function $w:\cU_f\to \R$ given by
 \[
 \U\mapsto \tr f(A(\U))
 \]
 is supermodular. \qed
 \end{theo}

 A function $w:\cU_f\to \R$ is {\em extendable} to $\cU$ if the following condition holds for each closed infinite dimensional subspace $\U\subset\cH$.
Let $\U_i, i\geq\1$
be an increasing sequence of finite dimensional subspaces of $\U$ such that $\cup_{i=1}^{\infty}\U_i$
 is dense in $\U$.
 Then the sequence $w(\U_i), i\geq 1$ converges to a unique value independent of the sequence $\U_i,i\geq 1$. We denote this limit by $w(\U)$.

 In the rest of this section, we discuss the extension of certain submodular and supermodular functions of a subspace $\cU$, arising from Theorem~\ref{theo-newmainid}.

 Denote by $\rC\rS_+(\cH)\subset \rS_+(\cH)$
 the closed ideal of positive compact operators.  For $A\in \rC\rS_+(\cH)$,
the Hilbert space $\cH$ has a countable orthonormal basis $\uu_1,\uu_2,\ldots$ such that
 \[A\uu_i=\lambda_i(A)\uu_i, i\geq 1, \quad \lambda_1(A)\ge \lambda_2(A) \ge \dots\ge 0, \lim_{i\to\infty}\lambda_i(A)=0.\]
 We now state some known results about characterizations
 of eigenvalues of $A\in\rS_+(\cH)$ that we will use here.  Consult with \cite{Fri73}.
 Let $\cX_m$ be the set of all orthonormal system of $m$ vectors  $X:=\{\x_1,\ldots,\x_m\}$ in $\cH$.
 The {\em convoy principle} (see Lemma~1 in~\cite{Fri73}) can be stated in the form of Corollary~4.4.3 in~\cite{Fbk}.
 \begin{equation}\label{convprin}
 \sup_{X\in\cX_m}\lambda_i(A(X))= \lambda_i(A) \textrm{ for } i=1,\ldots,m.
 \end{equation}
 The supremum is achieved for $X=\{\uu_1,\dots,\uu_m\}$.
 This in particular implies the Ky-Fan inequalities
 \[\tr (A(X))=\sum_{i=1}^m \an{A\x_i,\x_i}\le \sum_{i=1}^m \lambda_i(A).\]
 Let $\U=\span(\x_1,\ldots,\x_m)$.  Then the eigenvalues of $A(X)$ are the eigenvalues of $A(\U)$.
 The following approximation result is well known and can be straightforwardly deduced from the arguments in \cite{Fri73}.
 \begin{lemma}\label{approxprop}  Let $\cH$ be an infinite dimensional separable Hilbert space.  Let $A\in\rC\rS_+(\cH)$ and $\U$ be
 an infinite dimensional closed subspace of $\cH$.  Assume that $\U_i,i\geq 1$ is an increasing sequence of finite dimensional subspaces
 in $\U$, such that $\cup_{i\geq 1} \U_i$ is dense in $\U$. Then for each $j\geq 1$, the sequence $\lambda_j(A(\U_i)), i\geq 1,\dim\U_i\ge j$ is nondecreasing and converges to $\lambda_j(A(\U))$.
 \end{lemma}

 Let $p>0$.  Denote by $w_p: \cU_f\to \R$ the function $w_p(\U):=\tr A(\U)^p$ for $A\in\rS_+(\cH)$.
 The interlacing properties of hermitian matrices yield that $w_p$ is a nondecreasing function,
 i.e. $w_p(\U)\le w_p(\V)$ for $\U\subset \V$.
 For $p>0$ denote by $\rC\rS_{+,p}(\cH)$ the subset of all positive compact operators $A$ on $\cH$ such that
 $\tr A^p:=\sum_{j=1}^{\infty} \lambda_j(A)^p<\infty$.  Clearly $\rC\rS_{+,p}(\cH)\subsetneqq \rC\rS_{+,q}(\cH)$ for $0<p <q$.
 \begin{lemma}\label{pextend}
 Let $\cH$ be an infinite dimensional separable Hilbert space.  Assume that $p>0$ and $A\in\rC\rS_{+,p}(\cH)$.
 Then the function $\U\mapsto \tr A(\U)^p$ is extendable from $\cU_f$ to $\cU$. Moreover $w_p$ is a nondecreasing function.
 \end{lemma}
 \proof
 Let $\U_i, i\geq 1$ be a sequence of finite dimensional subspaces of $\U$ such that $\cup_{i=1}^{\infty} \U_i$ is dense in $\U$.
 Assume that $\dim\U_i=n_i$.  Then the positivity of $A(\U_i)$ and \eqref{convprin} imply that $0\le\lambda_j(A(\U_i))\le \lambda_j(A(\U))\le \lambda_j(A)$,
for $j=1,\ldots,n_i$.
Consider now, for each $i\geq 1$, the function
$\varphi_i:\{1,2,\dots\}\to [0,\infty)$ such that
$\varphi_i(j)=\lambda_j(A(\U_i))^p$ for $1\leq j\leq n_i$
and $\varphi_i(j)=0$ for $j>n_i$. Consider also the function
$\varphi:\{1,2,\dots\}\to [0,\infty)$ such that
$\varphi(j):=\lambda_j(A(\U))^p$.
Then, by Lemma~\ref{approxprop}, the sequence of nonnegative functions $\varphi_i, i\geq 1$ converges monotonically to the function $\varphi$ as $i\to\infty$, and then, the monotone convergence theorem implies that $\tr A(\U_i)^p=\sum_{j\geq 1}\varphi_i(j)$ converges, as $i\to\infty$, to $\sum_{j\geq 1}\varphi(j)=\tr A(\U)^p \leq \tr A^p<\infty$.
Finally, the convoy principle implies that $w:\cU\to\R$ is a nondecreasing function.  \qed

 \begin{lemma}\label{pextendentropy}
 Let $\cH$ be an infinite dimensional separable Hilbert space.
Assume that $A\in\rC\rS_{+}(\cH)$, and that $-\tr(A\log A)<\infty$.
 Then the function $\tilde{w}:\U\mapsto -\tr (A(\U)\log A(\U))$ is extendable from $\cU_f$ to $\cU$.
 \end{lemma}
 \proof
Let $\U_i,i\geq 1$ be as in the proof of Lemma~\ref{pextend}, still
with $n_i=\dim \U_i$.
Observe first that the function $-x\log x$ is increasing and nonnegative
on $[0,e^{-1}]$.
As $\lim_{j\to\infty} \lambda_j(A)=0$, we can find an index $\bar{\jmath}$ such that $\lambda_j(A)\leq e^{-1}$ for all $j\geq \bar{\jmath}$.
Let us fix such a $j$.
Then, by the convoy principle,
we have $0\leq \lambda_j(A(\U_i))\leq \lambda_j(A)\leq e^{-1}$, for all $i$,
and since the sequence $\lambda_j(A(\U_i)),i\geq 1$ is nondecreasing, it follows
that the sequence $ -\lambda_j(A(\U_i))\log\lambda_j(A(\U_i)), i\geq 1$
is nondecreasing and nonnegative. Applying the monotone convergence theorem
as in the proof of Lemma~\ref{pextend}, we deduce
that
\begin{eqnarray}
0 \leq \lim_{i\to\infty} \sum_{j\geq \bar{\jmath}} -\lambda_j(A(\U_i))\log\lambda_j(A(\U_i))
&=& \sum_{j\geq \bar{\jmath}} -\lambda_j(A(\U))\log\lambda_j(A(\U))\nonumber\\
&\leq& \sum_{j\geq \bar{\jmath}} -\lambda_j(A)\log\lambda_j(A) \enspace ,
\label{e-boundentropy}
\end{eqnarray}
again with the convention that $\lambda_j(\U_i)=0$ for $j>n_i$.
Since the sum in~\eqref{e-boundentropy} differs
from the sum
\( \sum_{j\geq 1}-\lambda_j(A)\log\lambda_j(A)=-\tr(A\log(A)) <\infty
\)
only by a finite number
of terms, we conclude that the sum in~\eqref{e-boundentropy} cannot be equal to $\infty$.
Since for all $0\leq j<\bar{\jmath}$, we also have $\lambda_j(A(\U_i))\to
\lambda_j(A(\U))$ as $i\to\infty$, we finally get
\[
-\infty<\lim_{i\to\infty} -\tr (A(\U_i)\log A(\U_i)) =
-\tr (A(\U)\log A(\U)) <\infty\enspace .
\]
\qed
 \begin{theo}\label{coro-powerentropop}
 Let $\cH$ be an infinite dimensional separable Hilbert space.  Assume that $p>0$ and $A\in\rC\rS_{+,p}(\cH)$.
 Then the function $w_p:\cU\to \R$ given by $w_p(\U):=\tr (A(\U))^p$ is nondecreasing, it is submodular for $p\in (0,1)$, supermodular
 for $p\in (1,2]$ and modular for $p=1$.  Furthermore,
if $A\in\rC\rS_{+}(\cH)$ is such that
 $-\tr(A \log A)<\infty$,
%
and in particular, if $A\in\rC\rS_{+,p}(\cH)$
for some $p\in (0,1)$, then the function $\tilde{w}(\U):=-\tr (A(\U)\log A(\U)), \U\in\cU$ is submodular.
 \end{theo}
 \proof The results for the function $w_p(\cdot)$ follows from Corollary \ref{coro-powerentrop} and Lemma \ref{pextend}. The submodularity
property for the function $\tilde{w}(\cdot)$ follows from
the same corollary and Lemma~\ref{pextendentropy}. Finally,
for $p\in (0,1)$, $(-x\log x)/x^p\to 0$ as $x\searrow 0$,
and since $\lambda_j(A)\searrow 0$ as $j\nearrow\infty$,
the convergence of the series $\tr A^p=\sum_{j} \lambda_j(A)^p$
implies the convergence of the series $-\tr(A\log A)= \sum_{j} -\lambda_j(A)\log \lambda_j(A)$. Hence, $-\tr(A\log(A))<\infty$ holds
as soon as $A\in \rC\rS_{+,p}(\cH)$ for some $p\in (0,1)$.
\qed
 \section{CUR approximation of nonnegative definite matrices}
 Let $A\in\C^{m\times n}$ be of rank $r=\rank A$.  Assume that $k\in [r-1]$.  The best rank $k$ approximation of $A$, denoted as $A_k$, is given by the \emph{singular value decomposition}, abbreviated here as SVD, \cite{GVL96}.  When $m$ and $n$ are very big, e.g. $n,m\ge 10^6$, finding $A_k$ is not computationally possible, since the complexity of computing $A_k$ is $O(kmn)$.  A good alternative for $A_k$ is the $CUR$ approximation,
 where $C\in\R^{m\times k}, R\in \C^{k\times n}$ are submatrices of $A$ of the form $A[[m],J],A[I,[n]]$, where $J\subset [n],I\subset [m]$ are subsets of columns and rows of $A$ of cardinality $k$.  Here $U\in\C^{k\times k}$ is a suitably chosen matrix.
 The best choice $U^{\star}$ for the Frobenius norm, i.e. $\arg \min_U \tr((A-CUR)(A-CUR)^*)$, is given by $C^{\dagger}AR^{\dagger}$, where $B^{\dagger}$
 denotes the Moore-Penrose inverse \cite{FMMN}.  Again, for $m,n\gg 1$ it is unfeasible to compute $U^{\star}$.

 Assume that $A[I,J]\in \C^{k\times k}$, the submatrix of $A$ based on the rows and columns $I,J$ respectively, is invertible.
 Then $CA[I,J]^{-1}R$ has the same $I$ rows and $J$ columns as $A$.  Hence the best choice of $I^{\star},J^{\star}$ seems to be given
 by~\cite{GTZ97}
 \begin{equation}\label{maxdetIJ}
 \mu_k(A):=\max \{|\det A[I,J]|, I\subset [m], J\subset [n], |I|=|J|=k\}=|\det A[I^{\star},J^{\star}]|.
 \end{equation}
 More precisely, for $A=[a_{ij}]_{i=j=1}^{m,n}\in\C^{m\times n}$ let $\|A\|_{\infty,e}:=\max_{i\in[m],j\in[n]}|a_{ij}|$
 be the entrywise max-norm of $A$.  Denote by $\sigma_1(A)\ge\ldots\ge \sigma_r (A)>0$ the $r$-positive singular values
 of $A$.  If $\det A[I,J]\ne 0$ one has the inequality \cite{GTZ97}
 \begin{equation}\label{CURest}
 \|A-CA[I,J]^{-1}R\|_{\infty,e}\le \frac{\mu_k(A)}{|\det A[I,J]|}(k+1)\sigma_{k+1}(A).
 \end{equation}
 (See \cite[\S4.13]{Fbk} for a simple proof of the above inequality.)

 In this section we discuss the maximal problem \eqref{maxdetIJ} and the related problem of $CUR$ approximation where $A\in\rH_{m,+}$.
 \begin{lemma}\label{maxdetpd}  For $A\in\rH_{m,+}$ and $k\in [m]$, the maximum in \eqref{maxdetIJ} is achieved for $I^{\star}=J^{\star}$, i.e.
 \begin{equation}\label{maxdetpdh}
 \mu_k(A)=\max_{I\subset [m], |I|=k} \det A[I].
 \end{equation}
 \end{lemma}
 \proof  Let $A=[a_{ij}]_{i=j=1}^m$.  We first consider the easy case $k=1$.  Recall that $a_{ii}\ge 0$ and each $2\time 2$ principal minor
 of $A$ is nonnegative.  Hence $a_{ii}a_{ij}-|a_{ij}|^2\ge 0$.  Thus $\max(a_{ii},a_{jj})\ge |a_{ij}|=|a_{ji}|$.  Hence \eqref{maxdetpdh} holds for $k=1$.  Assume that $k>1$.  Let $C_k(A)=\wedge^k A\in \C^{{m \choose k}\times {m \choose k}}$ be the $k$-th compound of $A$.  (The entries of $C_k(A)$ are all $\det A[I,J]$, $|I|=|J|=k$ arranged in a lexicographical order.)
 Recall that $C_k(A)\in \rH_{{m \choose k},+}$.  Hence the previous argument
applied now to any $2\times 2$ principal minor of the matrix $C_k(A)$
 yields \eqref{maxdetpdh}.  \qed

 We now introduce a simple
greedy algorithm for finding an approximation of $\mu_k(A)$ for $A\in \rH_{m,+}$ and the corresponding row index $I$.
 \begin{algo}\label{greedalg} $\;$

 \emph{Input}: $A=[a_{ij}]\in \rH_{m,+}, k\in [m], I=\emptyset, K=[m], t=1 $

 \emph{Output}: A subset $I$ of $[m]$ of cardinality $k$ and $\det A[I](=t)$

 $\;$\emph{Find} $i\in K$ $a_{ii}=\max_{j\in K} a_{jj}$.

 $\;\;$ \emph{If} $a_{ii}=0$ set $I=[k], t=0$ \emph{exit}

 $\;\;$ $I=I\cup\{i\}, K=K\setminus\{i\}, t=t a_{ii}$

 $\;\;$ \emph{If} $|I|=k$ \emph{exit}

 $\;\;$ \emph{For} $p,q\in K$ $a_{pq}=a_{pq}- \frac{a_{pi}a_{iq}}{a_{ii}}$

 \end{algo}

To analyse this algorithm, assume for the simplicity of argument that $\det A>0$.
Then, we may introduce the submodular function $w(I):=\tr (\log A[I])$,
and consider, as in the introduction, the maximum $\nu_k(\submodularf)$
of the function $w$ over all subsets $J\subset [m]$ of cardinality $k$,
so that
\(\nu_k(\submodularf) = \log \mu_k (A)
.
\)

Algorithm~\ref{greedalg} should be compared with the standard greedy
algorithm~\cite{NWF78} to approximate the maximum $\nu_k(\submodularf)$,
which consists in constructing
a sequence of sets $\emptyset=I^0\subset I^1\subset \dots\subset I^k\subset [m]$ such that at each step $r$, $I^r=I^{r-1}\cup\{i_r\}$ where the index
$i_r$ is chosen to so that $w(I^{r-1}\cup\{i_r\})= \max_{i\in [m]\setminus I^{r-1}} w(I^{r-1}\cup\{i\})$. We denote by $\nu^G(\submodularf):=\submodularf(I^k)$
the value of the solution $I^k$ returned by this greedy algorithm.
Algorithm~\ref{greedalg} is nothing but an implementation of this greedy algorithm, in which
at each step, the augmenting index $i_r$ is obtained as a maximal pivot
(third line of the algorithm)
and the value $w(I^{r})$ is obtained incrementally from $w(I^{r-1})$ by one
Gaussian elimination step (fifth and seventh lines of the algorithm).
In particular, denoting by $\mu^G_k(A)$ the value returned by Algorithm~\ref{greedalg},
we get
\(
\nu_k^G(\submodularf)= \log \mu_k^G(A) 
.
\)

Assume now that $w$ is nondecreasing.
Then, it is straightforward to show that if $m\ge 2$,
then $\det A[J]\ge 1$ for each  $J\subset [m]$, so that $w$ is
nonnegative.
Then, the result of~\cite{NWF78} can be applied to $w$,
showing that the greedy algorithm has an approximation factor $1-e^{-1}$,
i.e.,
$\nu_k^G(\submodularf)\ge (1-e^{-1})\nu_k(\submodularf)$.
It follows that the set $I$ returned by Algorithm~\ref{greedalg} satisfies
 \begin{equation}\label{greedest}
 \det A[I]=\mu_k^G(A)\ge \mu_k(A)^{1-e^{-1}} \enspace .
 \end{equation}
Therefore, we arrive at the following estimate of the CUR approximation
error of the greedy algorithm.
\begin{prop}
Let $m\geq 2$, and let $A\in \rH_{m,+}$ be such that $w(I):=\tr\log A[I]$ is nondecreasing. Then, the subset $I$ returned by Algorithm~\ref{greedalg}
satisfies:
\[
 \|A-A[[m],I]A[I]^{-1}A[[m],I]^*\|_{\infty,e}\le(k+1)\mu_k(A)^{e^{-1}}\sigma_{k+1}(A).
\]
\end{prop}
\proof
 Combine Inequality~\eqref{greedest} with Inequality~\eqref{CURest}.
\qed

 We now give a simple condition for $w$ to be a nondecreasing function.
 \begin{lemma}\label{wcon}  Let $A\in\rH_{m,+}$.  Assume that $\lambda_m(A)\ge 1$.  Then $w$ is a nondecreasing function.
 \end{lemma}
 \proof
 Let $I$ be a nontrivial subset of $m$ and let $i\in [m]\setminus I$.  As the eigenvalues of $B:=A[I]$ interlace the eigenvalues of $C:=A[I\cup\{i\}]$
 we deduce that
 \[\det C=\prod_{j=1}^{|I|+1} \lambda_j(C)\ge (\det B)\lambda_{|I|+1}(C)\ge (\det B)\lambda_m(A)\ge \det B.\]
 \qed

 We conclude this section with an upper estimate of the $CUR$ approximation error of $A\in\rH_{m,+}$ using Algorithm \ref{greedalg}. This estimate relies on the Hadamard determinant inequality and remains valid even when $w$ is not nondecreasing.
 \begin{prop}\label{ubgral} Let $A\in\rH_{m,+}$ and $k\in [m-1]$.  Assume that $k\le \rank A$.  Let $I$ be any subset of $[m]$ of cardinality $k$ returned by Algorithm \ref{greedalg}.  Let $a_1,\ldots,a_k$ be the $k$ maximal entries out of the $m$ diagonal entries of $A$.
 Then
 \[\|A-A[[m],I]A[I]^{-1}A[[m],I]^*\|_{\infty,e}\le \frac{\prod_{i=1}^k a_i}{\det A[I]}(k+1)\sigma_{k+1}(A).\]
 \end{prop}
 \proof  Let $I^{\star}$ be such that $\mu_k(A)=\det A[I^{\star}]$.  The Hadamard determinant inequality yields that $\det A[I^{\star}]\le \prod_{j\in I^{\star}}
 a_{jj}$.  Clearly  $\prod_{j\in I^{\star}} a_{jj}\le a_1\dots a_k$.  Use~\eqref{CURest} to deduce the lemma.\qed

 \section{Examples and further comments}\label{sec-further}

 We now give counter-examples showing that several natural generalizations
of the present results do not hold. In particular, we shall see
that the assumptions in Corollary~\ref{coro-powerentrop}
under which the function $I\mapsto T_p(I,A):=\tr A[I]^p$ is sub or super-modular,
for a $m\times m$ nonnegative definite matrix $A$, are tight.

To this end, we shall consider a hermitian $m\times m$ matrix $B$ with a block partition of the form
 \[
 B=\bordermatrix{&L_1&L_2&L_3\cr
 L_1& B_{11}&B_{12}&B_{13}\cr
 L_2&B_{12}^*&B_{22}&B_{23}\cr
 L_3&B_{13}^*&B_{23}^*&B_{33}}
 \]
 so that $[m]$ is the disjoint union of the sets $L_1,L_2,L_3$.
 Given $I,J\subset [m]$,
 it is convenient to set
 $$\delta(f,I,J,B):=\tr f(B[I])+\tr f(B[J])-\tr f(B[I\cup J])-\tr f(B[I\cap J])
 \enspace .$$
 Hence, $\delta(f,I,J,B)\geq 0$ for all $I,J$ whenever $I\mapsto \tr f(B[I])$ is
 submodular.

 We now choose
 \[
 I:=L_1\cup L_2,\qquad J:=L_2\cup L_3 \enspace .
 \]
 Using the interpretation of the matrix product in terms of concatenation of paths, we readily get
 \begin{equation}
 \delta(t^2,I,J,B)=-2\tr (B_{13}B_{13}^*)\leq 0 \enspace .\label{e-p2}
 \end{equation}
 Moreover, the inequality is strict as soon as $B_{13}\neq 0$.
 Since $\delta(t^p,I,J,B)$ depends continuously of $p$,
 we get the following proposition, in which one can take
$A=B$.
 \begin{prop}\label{prop-pnear2}
 For all $m\geq 3$, and for $p$ close enough to $2$,
 the strict supermodularity inequality for $T_p(\cdot,A)$ may hold.
 \qed
 \end{prop}
 In order to give examples in which the strict opposite inequality may hold,
 we observe that
 \begin{equation}
 B_{13}=0\implies
 \delta(t^3,I,J,B)=0 \quad\text{and}\quad \delta(t^4,I,J,B)=-4\tr B_{12}^*B_{12}B_{23}B_{23}^*
 \enspace .\label{e-jacobi}
 \end{equation}
 Note that $\delta(t^4,I,I,B)<0$, unless the nonnegative definite
 matrix $B_{12}B_{23}B_{23}^*B_{12}^*$ is zero. We assume
 in the sequel that this is not the case.

 Let us consider
 \[
 A = \Id + s B \enspace .
 \]
 Using~\eqref{e-p2} and~\eqref{e-jacobi},
  together with the trivial fact that $\delta(t,I,J,B)=0$,
 and $\tr (I+sB)^p=\tr I + ps\tr B+ \dots + \frac{p(p-1)(p-2)(p-3)}{4!}s^4\tr B^4
 +O(s^5)$, we identify the first terms of the
 Taylor expansion of $\delta(t^p,I,J,A)$ in $s$,
 \begin{eqnarray}
 \delta(t^p,I,J,A)&= -\frac{\displaystyle p(p-1)(p-2)(p-3)}{\displaystyle 3!}s^4\tr  B_{12}B_{12}^*B_{23}B_{23}^* +O(s^5)
 \end{eqnarray}
 For $s$ sufficiently small, the sign of the latter expression
 is positive for $p\in (2,3)$, and negative for $p>3$ or $p<0$.
 \begin{prop}\label{prop-23}
 For all $m\geq 3$, and for all $p\in (2,3)$ (resp $p<0$ or $p>3$)
 the strict submodularity (resp.~supermodularity) inequality for $T_p(\cdot,A)$ may hold.  \qed
 \end{prop}
 \begin{rem}\rm\label{rk-bapatsunders}
 The absence of systematic inequalities, except for $p\in [0,1]$
 or $p\in [1,2]$ is related to the absence of majorization
 properties. Indeed, consider again the matrix $B$ above,
 assume that it is of size $3\time 3$, so that the blocks
 $B_{12}$ and $B_{23}$ are scalars, and assume that $B_{13}=0$,
 so that $B$ is a Jacobi matrix. We choose again $I=\{1,2\}, J=\{2,3\}$.
 The spectra of the matrices $B,B[I],B[J],B[I\cap J]$ are
 \[\big(\sqrt{B_{12}^2+B_{23}^2},0, -\sqrt{B_{12}^2+B_{23}^2}\big),\quad \big(|B_{12}|,-|B_{12}|\big), \quad \big(|B_{23}|,-|B_{23}|\big), \quad \big(0\big) \enspace.\]
 Consider $A:=\Id+B$, so that $A>0$ for $1> \sqrt{B_{12}^2+B_{23}^2}$.
  The joint sequences of the eigenvalues of $A$ with $A[I\cap J]$ and of $A[I]$ with $A[J]$ are
  $$\big(\sqrt{B_{12}^2+B_{23}^2}+1,1,1,-\sqrt{B_{12}^2+B_{23}^2}+1\big),\quad \big(|B_{12}|+1,|B_{23}|+1, -|B_{23}|+1,-|B_{12}|+1\big).$$
  Note that for $B_{12}B_{23}\ne 0$ none of the joint eigenvalue sequences majorizes the other one.
 This example should be compared with a theorem of
 Bapat and Sunder~\cite{bapatsunder}, concerning the special case in which $I\cap J=\emptyset$: then, the sequence obtained by concatenating the eigenvalues of $A[I]$ and $A[J]$ is majorized by the sequence of eigenvalues of $A[I\cup J]$,
 and so, the inequality~\eqref{powerentro12} for the map $I\mapsto \tr A[I]^p$
 holds for for all $p\in (-\infty,0) \cup [1,\infty)$ if $I\cap J=\emptyset$.
 \end{rem}

 \begin{example}\rm\label{ex--1}
 There are negative values of $p$ for which $T_p(\cdot,A)$ is not
 supermodular. Consider for instance the positive definite matrix
 \[
 A=\left(\begin{array}{rrr}
    5&   - 12&    9\\
   - 12&    33&  - 24\\
     9&   - 24&    19
 \end{array}
 \right),
 \quad I=\{1,2\},J=\{1,3\},\quad \delta(t^{-1},I,J,A)=16/35>0 \enspace .
 \]
 \end{example}

 The following counter example shows that Theorem~\ref{theo-newmain} does not carry over to operator convex maps.
 \begin{example}\rm
 The map $f(t)=t^2(\lambda+t)^{-1}$ is known to be operator convex~\cite[Problem~V.5.5]{Bha97}. Consider the positive definite matrix
 \[
 A=\left(\begin{array}{rrr}
     1&   - 2&  - 2\\
   - 2&    6&    4\\
   - 2&    4&    8
 \end{array}\right),
 \qquad \lambda=1,\qquad I=\{1,2\}, J=\{1,3\} \enspace .
 \]
 Then, $\delta(f,I,J,A)
 = 44/1085>0$, showing that the map $K\mapsto \tr f(A[K])$ is not supermodular.
 \end{example}

 We now summarize in Table~\ref{table1} the submodularity and supermodularity properties we know, omitting the proofs of the two easy cases.
 \begin{table}[htpb]
 \begin{center}
 \begin{tabular}{c|c|c}
 $p=-1$, $m\geq 3$& may not be supermodular & Ex.~\ref{ex--1}\\
 $p<0$, $m\geq 3$& may not be submodular & Prop.~\ref{prop-23}\\
 $p<0$, $m=2$& supermodular & Omitted\\
 $0\leq p\leq 1$& submodular & Coro.~\ref{coro-powerentrop}\\
 $1\leq p\leq 2$& supermodular& Coro.~\ref{coro-powerentrop}\\
 $p>2$, $m=2$& supermodular & ~Omitted\\
 $p=2$, $m\geq 3$ & may not be submodular & Prop.~\ref{prop-pnear2}\\
 $p\in (2,3)$, $m\geq 3$ & may not be supermodular &Prop.~\ref{prop-23}\\
 $p>3$, $m\geq 3$ & may not be submodular &Prop.~\ref{prop-23}
 \end{tabular}
 \end{center}
 \caption{Summary of submodularity and supermodularity properties of the spectral function $T_p(I,A)=\tr A[I]^p$ for a $m\times m$ nonnegative definite Hermitean matrix $A$.}\label{table1}
 \end{table}

 Recall that for $A,B\in \rH_m$ and $A\le B$ we have that $\tr A\le \tr B$.
 Observe also that $\rH_m(\interval)$ is a convex subset of $\rH_m$. Since the primitive of an operator monotone function is operator convex we deduce from Theorem~\ref{theo-newmain}:
 \begin{corol}\label{trconv}  Let $f$ satisfies the assumptions of Theorem \ref{theo-newmain}.
 Then for each $I\subseteq [m]$ the function $A\mapsto \tr f(A[I])$ is a convex function on
 $\rH_m(\interval)$.\qed
 \end{corol}

 \begin{rem}\rm
 A different submodularity inequality involving spectral functions
 appeared in~\cite{rio}, with an application to an experiment design problem.
 It is shown there that for all nonnegative definite
 hermitian matrices $A,B,C$, and for all $0\leq p\leq 1$,
 $\tr (A+B+C)^p+\tr C^p\leq \tr (A+C)^p + \tr (B+C)^p$.
 \end{rem}

 \end{document}